\magnification=1200
{\centerline{\bf Uniqueness of Solutions to Schr\"odinger Equations on
Complex Semi-Simple Lie Groups}}
\bigskip
{\centerline{\bf  by}}
\bigskip
{\centerline{\bf Sagun Chanillo}}
\bigskip
\bigskip
\bigskip
{\centerline{\bf Abstract}}
\bigskip
In this note we study the time dependent Schr\"odinger equation on
Complex semi-simple Lie Groups. We show that if the initial data is a
bi-invariant function that has
sufficient decay and the solution has sufficient decay at another
fixed value of time, then the solution has to be identically zero for
all time. We
also derive Strichartz and decay estimates for the Schr\"odinger
equation. Our methods also extend to the wave equation. On the
Heisenberg group we show that the failure to obtain a parametrix for
our Schr\"odinger equation is related to the fact that geodesics
project to circles on the contact plane at the identity.
\bigskip
\medskip
{\centerline{\it Dedicated to U. B. Tewari on his Retirement}}
\bigskip
\bigskip
\bigskip
\bigskip
\noindent
MSC(Subj. Classification): 43A85
\bigskip
\bigskip
\noindent
Running Head: Schr\"odinger's Equation on Complex Lie Groups
\vfill\eject
\noindent
{\bf \S 1.  Introduction:} 

Let $G$ denote a Lie group. We are concerned here with
the initial value problem for the time dependent Schr\"odinger eqn.,
    $$-iu_t=\Delta u,\ u|_{t=0}=f(x),\ x\in G\eqno (1.1)$$
$\Delta$ will denote the
Laplace-Beltrami operator with respect to an appropriate $G$ invariant
metric. Our main aim is to obtain 
conditions on $f$ which will guarantee that the solution $u(x,t)$,
vanishes identically for all $t>0$. The results that we have may be
all viewed as consequences of a well-known theorem of Hardy on the
Euclidean Fourier transform, which is a form of the uncertainty
principle. The Hardy theorem has been extended to the Lie group
setting by various authors, see [2], [9]. However, we do not need
any Lie group version of the Hardy theorem, the original Euclidean
version will suffice in our applications to Lie groups. We plan to extend our results to the real
semi-simple Lie groups and to obtain versions on the Heisenberg
group later. On nilpotent groups there is a serious difficulty since the
geodesics project to circles on the contact plane and the existence of
such closed loops creates difficulties in writing down a parametrix
for the solution operator to $(1.1)$, see section 3. Our results may
be viewed as a statement that a concentrated wave packet at initial
time will spread out and if still remains concentrated, it must be
trivial. Lastly in Sec 4. we derive decay estimates and Strichartz
estimates for bi-invariant solutions to the Schr\"odinger equation on
complex semi-simple groups. The methods extend to the wave equation on
complex Lie groups.

We wish to thank Roe Goodman for reading our paper and offering
constructive criticism. The research for this paper was supported in
part by an NSF grant DMS-0600971.
\bigskip
\noindent
{\bf \S 2. The Uniqueness theorem}

We recall the theorem of Hardy. Let the Fourier transform be defined
as,
  $$\hat f(\xi)=\int_{R^n}e^{-i<x,\xi>}f(x)\, dx.$$
\medskip
\noindent
{\bf Theorem 1:}(Hardy [4]) Let $f(x)$ satisfy $|f(x)|\leq
Ce^{-a|x|^2}$. Furthermore,
assume, $|\hat f(\xi)|\leq Be^{-b|\xi|^2}$. If $4ab>1$, then $f\equiv
0$.
\bigskip
We shall now derive from this theorem a uniqueness theorem for
solutions to $(1.1)$ when $G=R^n$. 
\medskip
\noindent
{\bf Theorem 2:} Let us consider the initial value problem for $u(x,t)$,
    $$-iu_t=\Delta u,\  u|_{t=0}=f(x).$$

Assume that $|f(x)|\leq Ae^{-a|x|^2}$, and $|u(x,t_0)|\leq B
e^{-b|x|^2}$. If $16abt_0^2>1$, then $u(x,t)\equiv 0$, for all $t\geq 0$.

\medskip
\noindent
{\bf Proof:} We recall that using the fundamental solution of $(1.1)$,
we may write,
    $$u(x,t)={{c_n}\over{t^{n/2}}}\int_{R^n}e^{i{{|x-y|^2}\over{4t}}}f(y)\, dy.$$
We may re-write the last identity as,
  $$u(x,t)={{c_ne^{i{{|x|^2}\over{4t}}}}\over{t^{n/2}}}\int_{R^n}
e^{-i<{{x}\over{2t}},y>}e^{i{{|y|^2}\over{4t}}}f(y)\, dy.\eqno (2.1)$$
Now set $h(y)=e^{i{{|y|^2}\over{4t}}}f(y)$. Then from $(2.1)$ we get,
    $$u(x,t)= {{c_ne^{i{{|x|^2}\over{4t}}}}\over{t^{n/2}}}
\hat h({{x}\over{2t}}).$$
Now we take $t=t_0$ and apply Hardy's theorem. From the hypothesis on
$u(x,t_0)$, we have,
    $$|\hat h ({{x}\over{2t_0}})|\leq Be^{-b|x|^2}$$
Thus,
      $$|\hat h(x)|\leq Be^{-4bt_0^2|x|^2}$$ 
Clearly $|h(x)|\leq C e^{-a|x|^2}$. Since $16abt_0^2>1$, Hardy's
theorem applies and we
conclude $h\equiv 0$. It then follows that $f\equiv 0$ and hence
$u(\cdot , t)\equiv 0$ for all $t$.

\bigskip
We now wish to extend Theorem 2 to a complex semi-simple Lie
group. We need to introduce some notation. 

Let $G$ denote a complex, connected  semi-simple Lie group. Let $K$ denote a fixed maximal compact sub-group of $G$. Let
${\cal G}$ and ${\cal K}$ denote the Lie algebras of $G$ and $K$
respectively. Let $B$ denote the Cartan-Killing form on ${\cal G}$,
and the Cartan decomposition is given by ${\cal G}={\cal K}
\oplus{\cal P}$. The restriction of $B$ to ${\cal P}\times {\cal P}$ is strictly positive definite and hence defines a norm. Let ${\cal A}$
be a fixed maximal Abelian subspace of ${\cal P}$. Let $\Sigma$
denote the set of non-zero roots corresponding to the pair $({\cal
G},{\cal A})$, and $\Sigma_+$ the set of positive roots $\alpha$ for some
ordering. $W$ will denote the Weyl group associated to $\Sigma$. Let
${\cal A}_+$ be the positive Weyl chamber associated to $\Sigma$. Let
$A$ denote the analytic sub-group with Lie algebra ${\cal A}$. We
have the map, $\exp:\ {\cal A}\to A$. Likewise we have,
$A_+=\exp({\cal  A}_+)$. Then one has the polar decomposition,
   $$G=K\overline{A_+}K\eqno (2.2)$$
The Haar measure on $A$ can then be written for this polar
decomposition by a formula of Harish-Chandra [5] as,
  $$dx=(\sum_{s \in W }( \det s)e^{s\rho (H)})^2\, dH=\phi(H)^2\, dH, \ H\in {\cal
A}.\eqno (2.3)$$
As usual $\rho={{1}\over{2}}\sum_{\Sigma_+}\alpha$. The spectral
variables are elements of ${\cal A}^\star$ and will be denoted by
$\lambda$. On $G/K$ we have a $G$ invariant Riemannian metric obtained through
the Killing form and we can form a Laplace-Beltrami operator $\Delta$ using
this metric. See [7]. Note now, that there is a unique element
$H_\lambda\in {\cal A}$, such that,
    $$\lambda(H)=B(H_\lambda, H), \ H\in {\cal A}.\eqno (2.4)$$
  The norm of elements $H\in {\cal A}$ will be denoted by,
      $$|H|^2=B(H,H)$$    
Lastly a function $f$ on $G$ will be said to be $K$ bi-invariant if
and only if,
         $$f(k_1ak_2)=f(a),\ a \in A$$
From the polar decomposition we may view the function as only
depending on its values on $A_+$, or by using the inverse exponential
map we may also view $f$ as a complex valued function on ${\cal A}$
solely determined by its values on ${\cal A}_+$. We are ready to
state our theorem.

\bigskip
\noindent
{\bf Theorem 3:} Let $G$ denote a connected, complex, semi-simple Lie
group. Let $f$ be a bi-invariant function such
that, 
  $$|f(H)|\leq A e^{-a|H|^2}.$$
Consider the initial value problem,
    $$-iu_t=\Delta u,\ u|_{t=0}=f$$
Then $u$ is also bi-invariant. Furthermore if
      $$|u(H,t_0)|\leq B e^{-b|H|^2},$$
and $16abt_0^2>1$, then necessarily $u\equiv 0$ for all $t\geq 0$.

\bigskip
\noindent
{\bf Proof:} The proof requires the use of well-known facts about
spherical functions. First we need to recall Theorem 5.7, p.432, [7]
which states that on a complex semi-simple group, the elementary
spherical functions are given by(where $\phi(H)$ is as in $(2.3)$),
    $$\phi_\lambda(H)= c(\lambda) {{\sum_{s\in W}(\det s)e^{is\lambda
(H)}}\over{\phi (H)}}\eqno (2.5)$$
Here $c(\lambda)$ is the celebrated c-function of Harish-Chandra [5,6], [3]. Furthermore we also have from [7] that,
  $$\Delta \phi_\lambda(H)=-(|\lambda|^2+|\rho|^2)\phi_\lambda(H)\eqno (2.6)$$
The spherical transform of $f$ is given by,
  $$\hat f(\lambda)=\int_{{\cal A}}\phi_{-\lambda}(H)f(H)\phi^2(H)
dH. \eqno (2.7)$$

Thus the solution of our initial value problem in view of $(2.6)$ is,
    $$u(H,t)=\int_{{\cal
A}^\star}e^{-it(|\lambda|^2+|\rho|^2)}\phi_\lambda(H)\hat f(\lambda)
|c(\lambda)|^{-2}\, d\lambda.\eqno (2.8)$$
We insert $(2.5)$, $(2.7)$ into  $(2.8)$ to get,
  $$u(H_1,t)\phi(H_1)=\int_{\cal A}(\int_{{\cal A}^\star}
e^{-it(|\lambda|^2+|\rho|^2)}(\sum_{s\in W}(\det
s)e^{is\lambda(H_1)})(\sum_{s^\prime\in W}(\det s^\prime)e^{-is^\prime\lambda(H_2)})\,
d\lambda)$$
$$f(H_2)\phi(H_2)\, dH_2.\eqno (2.9)$$
The inner integral in $\lambda$ may be re-written in view of $(2.4)$
as,

    $$ \sum_{s,s^\prime}(\det s)(\det
s^\prime)e^{-it|\rho|^2}e^{i{{|sH_1-s^\prime H_2|^2}\over{4t}}}\int_{{\cal
A}}e^{-it|H_\lambda+{{1}\over{2t}}(sH_1-s^\prime H_2)|^2}\,
dH_\lambda.\eqno (2.10)$$
Evaluating the integral in $(2.10)$ we get from $(2.9)$, ($l=\dim
{\cal A}$), 
  $$u(H_1,t)\phi(H_1)= c_l{{e^{-it|\rho|^2}}\over{t^{l/2}}}\int_{\cal
A}(\sum_{s,s^\prime}(\det s)(\det s^\prime)e^{i{{|sH_1-s^\prime
H_2|^2}\over{4t}}})f(H_2)\phi(H_2)dH_2.\eqno (2.11)$$
Set $g(H_2)=f(H_2)\phi(H_2)$. Now note that because $f$ is
bi-invariant and $\phi(H)$ is odd under the action of the Weyl group,
$g(sH)=(\det s) g(H)$. We proceed to re-write $(2.11)$ as in
the Euclidean case and we will make use of the Euclidean Fourier
transform. Re-writing $(2.11)$ we get,
$$u(H_1,t)\phi(H_1)=c_l{{e^{-it|\rho|^2}e^{i{{|H_1|^2}\over{4t}}}}\over{t^{l/2}}}
\sum_{s,s^\prime}(\det
s)(\det s^\prime)\int_{R^l}
e^{-iB({{ss^\prime
H_1}\over{2t}},H_2)}e^{i{{|H_2|^2}\over{4t}}}g(H_2)dH_2.\eqno (2.12)$$
Here $l$ is the rank of $G$. Set

$$R(H_2)=e^{i{{|H_2|^2}\over{4t}}}g(H_2).\eqno (2.13)$$ 
Then $(2.12)$ states,
  $$u(H_1,t)\phi(H_1)=c_l{{e^{-i(t|\rho|^2-{{|H_1|^2}\over{4t}})}}\over{t^{l/2}}}\sum_{s,s^\prime}(\det s)(\det s^\prime)\hat R({{ss^\prime H_1}\over{2t}})\eqno (2.14)$$

Since $R(H)$ is odd under reflection by the Weyl group, we finally can
write $(2.14)$ as,
  $$u(H_1,t)\phi(H_1)=c_l|W|^2{{e^{-i(t|\rho|^2-{{|H_1|^2}\over{4t}})}}\over{t^{l/2}}}\hat R({{H_1}\over{2t}}).\eqno (2.15)$$

Noting further from $(2.3)$ that $|\phi(H)|\leq ce^{c|H|}$, we see
right away from $(2.15)$ and the hypothesis that at $t=t_0$, 
  $$|\hat R(H_1)|\leq c e^{-4b^\prime t_0^2|H_1|^2}, b^\prime<b$$
and also,
     $$|R(H_2)|\leq ce^{-a^\prime |H_2|^2}, a^\prime<a.$$
Thus by Theorem 1, we again conclude $R\equiv 0$. This implies
$f\equiv 0$ and hence the theorem follows. 
\bigskip
\noindent
{\bf Lemma 1:} The results of Theorems 2 and 3 are sharp. One cannot
relax the condition $16abt_0^2>1$.
\medskip
\noindent
{\bf Proof:} We only display the proof for $R^n$. For theorem 3 one
can check the validity of the lemma by doing an explicit computation
on $SL(2,C)$. For $R^n$, choose initial data, 
$$f(x)= e^{-|x|^2-i{{|x|^2}\over{4}}}$$
An elementary computation using $(2.1)$ shows that,
   $$u(x,1)= c_ne^{{{i|x|^2}\over{4}}}e^{-{{|x|^2}\over{16}}}$$
Thus we have $16abt_0^2=1$ and uniqueness fails.
\bigskip\bigskip
\noindent
{\bf \S 3. The Heisenberg Group.}
\bigskip
We may ask the uniqueness question above for the sub-Laplacian
$\Delta_b$ on the
Heisenberg group. That is we consider the Schr\"odinger equation,
$$iu_t=\Delta_b u, \ u|_{t=0}=f(x)\eqno (3.1)$$
However there is a difficulty in writing the fundamental solution to
this operator due to the presence of closed loops in the contact
plane. To see this we proceed heuristically. First of all the
fundamental solution to the heat equation on the Heisenberg group is
given by,
    $$K(x,u,\xi,t)=\int_{-\infty}^{\infty}e^{-i\lambda
  \xi}e^{-t\lambda^2}{{\lambda}\over{\sinh
    (\lambda t)}}e^{-{{1}\over{4}}\lambda \coth (\lambda t)(x^2+u^2)}\,
d\lambda\eqno (3.2)$$
where the Heisenberg group is viewed as $R^3$ and points on it written
as $(x,u,\xi)$. To try to get a solution operator for $(3.1)$, we perform a change of variables in $(3.2)$ by
letting, $t\to -it$. The integrand in $(3.2)$ becomes,
   $$
{{\lambda}\over{\sin
    (\lambda t)}}e^{-{{i}\over{4}}\lambda\cot (\lambda
    t)(x^2+u^2)}\eqno (3.3)$$
Thus we note that the putative solution operator for fixed $t$ is
    singular at $\lambda=k\pi/t$. In fact $(3.3)$ converges
    to a Dirac delta at $u=x=0$ as $\lambda\to k\pi/t.$ This phenomena
    is attributable to the geodesics of the Heisenberg group
    projecting onto circles(closed loops) in the contact plane at the
    origin and having a cut-locus at $k\pi/t$. The geodesics on the Heisenberg group are well-known and a
    formula is found in [8]. They are given by,

   $$x(s)= {{\cos \beta(1-\cos(t s))+\sin\beta\sin(t
s)}\over{t}}$$
    $$u(s)={{-\sin \beta(1-\cos (ts))+\cos\beta\sin (t
s)}\over{t}}$$
$$\xi(s)=2{{ts-\sin(ts)}\over{t^2}}$$
Now notice the geodesics are not closed but their projection into the
$x-u$ plane, the contact plane,  are circles given by,
$$(x(s)-{{\cos\beta}\over{t}})^2+(u(s)+{{\sin\beta}\over{t}})^2={{1}\over{t^2}}.$$
The cut-locus emerges exactly at a distance of $k\pi/t$ from $u=x=0$
along the circle of radius $1/t$, and it is exactly there the
integrand of the solutions operator has its singularities.
A general construction of the heat kernel on CR manifolds in [1]
exhibits the same phenomena.
\bigskip\bigskip\noindent
{\bf \S 4. Strichartz Estimates and Decay estimates on Complex Lie
Groups}
\bigskip
We will use the results of our computation in Sec. 2 on  complex
semi-simple groups to obtain various estimates on the solution
$u(H,t)$. The estimates fall into two categories. One where we
integrate over a fixed time slice and another where we integrate over
both space and time, the Strichartz estimates. We have,
\bigskip
\noindent
{\bf Theorem 4:} Under the assumptions of theorem 3, the solution
$u(H,t)$ for bi-invariant initial data $f$ satisfies,
$$ ||u|\phi|^{1-{{2}\over{q}}}||_{L^q(G)}\leq
ct^{-l({{1}\over{p}}-{{1}\over{2}})}||f|\phi|^{1-{{2}\over{p}}}||_{L^p(G)},\
1\leq p\leq 2,\  {{1}\over{p}}+{{1}\over{q}}=1.\eqno (a)$$
Assume that,
          $$-iu_t-\Delta u= \psi(H,t), u|_{t=0}=f,$$
where $\psi$ is also bi-invariant for every fixed $t$. Then for
$p=2(l+2)/(l+4),\ q=2(l+2)/l$, we have,
    $$||u|\phi|^{1-{{2}\over{q}}}||_{L^q(G\times (0,\infty))}\leq
c(||f||_{L^2(G)}+||\psi|\phi|^{1-{{2}\over{p}}}||_{L^p(G\times
(0,\infty))})\eqno (b).$$
\bigskip
\noindent
{\bf Proof:} The proofs follow from the identity $(2.11)$. It
follows from $(2.11)$ that the function $u(H,t)\phi(H)$ is obtained by
applying the Euclidean fundamental solution of the Schr\"odinger
equation to the data given by $g(H)$. Thus from the Euclidean
estimates,
              $$||u\phi||_{L^q(R^l)}
\leq ct^{-l({{1}\over{p}}-{{1}\over{2}})}||g||_{L^p(R^l)},$$
where $p,q$ is as in part $(a)$ above. Re-writing the last inequality
using $(2.3)$ we get part $(a)$ of our theorem. 

A similar computation as in Cor. 1, in [10] gives part $(b)$. 
\medskip
\noindent
{\bf Remark:} The methods of Theorem 4 and [10] also extend to the
wave equation on complex semi-simple Lie groups with bi-invariant
functions as data. Strichartz and decay
estimates are easily obtained by following the methods developed in
section 2 and theorem 4. 
\vfill\eject
{\centerline{\bf References}}
\bigskip
\noindent
[1] Beals, R., Greiner, P. C. and Stanton, N.,  The Heat Equation on a CR
Manifold, J. Differential Geom., {\bf 20}, (1984), 343-387.
\medskip
\noindent
[2] Folland, G. B.,  and Sitaram, A.,  The uncertainty Principles: A
Mathematical Survey, J. Fourier Analysis Appl. {\bf 3}, (1997),
207-238. 
\medskip
\noindent
[3] Gindikin, S.  and Karpelevic, F. I.,  Plancherel Measure of
Riemannian Symmetric spaces of non-positive curvature,
Dokl. Akad. Nauk. SSSR, {\bf 145}, (1962), 252-255.

\medskip
\noindent
[4] Hardy, G. H.,  A theorem concerning Fourier transforms, J. London
Math. Soc., {\bf 8}, (1933), 227-231.
\medskip
\noindent
[5] Harish-Chandra, Spherical Functions on a Semi-simple Lie Group
I, Amer. J. of Math., {\bf 80}, (1958), 241-310.

\medskip
\noindent
[6] Harish-Chandra, Spherical functions on a Semi-simple Lie group
II, Amer. J. of Math., {\bf 80}, (1958), 553-613.
\medskip
\noindent
[7] Helgason, S.,  Groups and Geometric Analysis-Integral Geometry,
Invariant Differential Operators and Spherical Functions, Academic
Press, New York, 1984.
\medskip
\noindent
[8] Monti, R.,  Some properties of Carnot-Caratheodory balls in the
Heisenberg group, Rend. Math. Acc. Lincei, {\bf 11}, s.9, (2000),
155-167.
\medskip
\noindent
[9] Sitaram, A., and Sundari, M., An analogue of Hardy's theorem for
very rapidly decreasing functions, Pacific J. Math., {\bf 177},
(1997), 187-200.
\medskip
\noindent
[10] Strichartz, R. S.,  Restrictions of Fourier Transforms to Quadratic
Surfaces and Decay of Solutions to Wave Equations, Duke Math. J. {\bf
44}(3), (1977), 705-714.
\bigskip
\bigskip
\noindent
Deptt. of Math.,

\noindent
Rutgers University,

\noindent
110 Frelinghuysen Rd.,

\noindent
Piscataway, NJ 08854

\noindent
USA

\noindent
e-mail address: {\tt chanillo@math.rutgers.edu}
\end